\documentclass{elsart}
\usepackage{amsmath}
\usepackage{amssymb}

\newcommand{\R}{\mathbb{R}}                     

\begin{document}

\begin{frontmatter}

\title{Multivariate integration of functions depending explicitly on the minimum and the maximum of the variables}

\author{Jean-Luc Marichal}
\ead{jean-luc.marichal[at]uni.lu}

\address{
Institute of Mathematics, University of Luxembourg\\
162A, avenue de la Fa\"{\i}encerie, L-1511 Luxembourg, Luxembourg }

\date{Revised, September 8, 2007}

\begin{abstract}
By using some basic calculus of multiple integration, we provide an alternative expression of the integral
$$
\int_{]a,b[^n}f(\mathbf{x},\min x_i,\max x_i) \, d\mathbf{x},
$$
in which the minimum and the maximum are replaced with two single variables. We demonstrate the usefulness of
that expression in the computation of orness and andness average values of certain aggregation functions. By
generalizing our result to Riemann-Stieltjes integrals, we also provide a method for the calculation of certain
expected values and distribution functions.
\end{abstract}

\end{frontmatter}

\section{Introduction}

Let $a,b\in\mathbb{R}\cup\{-\infty,+\infty\}$, with $a<b$, and consider an integral over $]a,b[^n$ whose
integrand displays an explicit dependence on the minimum and/or the maximum of the variables, that is, an
integral of the form
\begin{equation}\label{eq:MainInt}
\int_{]a,b[^n}f(\mathbf{x},\min x_i,\max x_i) \, d\mathbf{x}.
\end{equation}
In this note we provide an alternative expression of this integral, in which the minimum and the maximum are
replaced with two single variables. When the integral is tractable, that alternative expression generally makes
the integral much easier to evaluate. For instance, when the integrand depends only on the minimum and the
maximum of the variables, we obtain the following identity
\begin{equation}\label{eq:nMaxMin}
\int_{\left]a,b\right[^n}f(\min x_i,\max x_i) \, d\mathbf{x} = n(n-1) \int_a^b dv\int_a^v f(u,v) (v-u)^{n-2}\,
du,
\end{equation}
and hence, for certain functions $f$, the integral becomes very easy to evaluate.

The alternative expression we present for integral (\ref{eq:MainInt}) is given in the next section (see
Theorem~\ref{thm:main}). The method we employ to obtain that expression merely consists in dividing the domain
$]a,b[^n$ into $n$ polyhedra chosen in such a way that the minimum and maximum functions simply become single
variables.

This method can be very efficient in the evaluation of many integrals that would normally require difficult and
tedious computations. As an example, consider the variance of a sample $\mathbf{x}\in [a,b]^n$ from a given
population, namely
$$
s^2(\mathbf{x})=\frac 1{n-1}\sum_{i=1}^n \Big(x_i-\frac 1n\sum_{j=1}^nx_j\Big)^2.
$$
The average value over $[a,b]^n$ of the variance-to-range ratio function can be easily calculated by using our
method. We merely obtain
\begin{equation}\label{eq:VarRan}
\frac 1{(b-a)^n}\int_{[a,b]^n}\frac{s^2(\mathbf{x})}{\max x_i-\min x_i} \, d\mathbf{x} = \frac{n+2}{12n}\,
(b-a).
\end{equation}

This note is set out as follows. In Section 2 we state and prove the main result. In Section 3 we provide an
application of our result to internal functions, also called Cauchy means, which can be classified according to
their location within the range of the variables. A similar application to conjunctive and disjunctive functions
is also investigated. In Section 4 we show how the direct generalization of our result to Riemann-Stieltjes
integrals enables us to consider the evaluation of certain expected values from various distributions.

We will use the following notation throughout. For any $n$-tuple $\mathbf{x}$, we denote by $(\mathbf{x}\mid
x_j=u)$ the $n$-tuple whose $i$th coordinate is $u$ if $i=j$, and $x_i$ otherwise. Also, for any integer
$n\geqslant 1$, we set $[n]:=\{1,\ldots,n\}$.

\section{Main result}

In this section we present our main result which consists of an alternative expression of integral
(\ref{eq:MainInt}). We start with a preliminary lemma, which concerns the particular cases of functions
involving either the minimum or the maximum of the variables.

\begin{lem}\label{lemma:main}
Let $f:\left]a,b\right[^{n+1}\to\mathbb{R}$ be an integrable function. Then we have
\begin{eqnarray*}
\int_{\left]a,b\right[^n}f(\mathbf{x},\min x_i) \, d\mathbf{x} &=& %
\sum_{j=1}^n \int_a^b du \int_{]u,b[^{n-1}} f(\mathbf{x},u\mid x_j=u) \, \prod_{i\in [n]\setminus \{j\}} dx_i,\\
\int_{\left]a,b\right[^n}f(\mathbf{x},\max x_i) \, d\mathbf{x} &=& %
\sum_{j=1}^n \int_a^b dv \int_{]a,v[^{n-1}} f(\mathbf{x},v\mid x_j=v) \, \prod_{i\in [n]\setminus \{j\}} dx_i.
\end{eqnarray*}
\end{lem}

\begin{pf*}{Proof.}
Consider the following $n$-dimensional open polyhedra
$$
P_j:=\{\mathbf{x}\in\,]a,b[^n\, : \, x_i>x_j\, \forall i\neq j\}\qquad (j\in [n]).
$$
They are pairwise disjoint. Indeed, if $\mathbf{x}\in P_j\cap P_k$, with $j\neq k$, then $x_k>x_j$ and
$x_j>x_k$, which is a contradiction. Moreover, the union of their set closures covers $]a,b[^n$. Indeed, for any
$\mathbf{x}\in\,]a,b[^n$ there is always $j\in [n]$ such that $x_i\geqslant x_j$ for all $i\neq j$.

Therefore, for any integrable function $f:\left]a,b\right[^{n+1}\to\mathbb{R}$, we have
\begin{eqnarray*}
\int_{\left]a,b\right[^n}f(\mathbf{x},\min x_i) \, d\mathbf{x} &=& \sum_{j=1}^n \int_{P_j}f(\mathbf{x},\min x_i) \, d\mathbf{x}\\
&=& \sum_{j=1}^n \int_a^b dx_j \int_{]x_j,b[^{n-1}} f(\mathbf{x},x_j) \, \prod_{i\in [n]\setminus \{j\}} dx_i,
\end{eqnarray*}
which proves the first formula. The second formula can be established similarly by considering the polyhedra
$$
Q_j:=\{\mathbf{x}\in\,]a,b[^n\, : \, x_i<x_j\, \forall i\neq j\}\qquad (j\in [n]).\qed
$$
\end{pf*}

Lemma~\ref{lemma:main} is interesting in its own right since it provides special cases of the main result. For
instance, by applying the first formula, we immediately obtain the following identity, which will be used in the
next section (see Example~\ref{ex:Choquet}). For any $S\subseteq [n]$, we have
\begin{equation}\label{eq:CaseMinSab}
\int_{]a,b[^n} f\big(\min_{i\in S}x_i\big)\, d\mathbf{x} =  (b-a)^{n-|S|}\, |S| \int_a^b f(u) (b-u)^{|S|-1}\,
du.
\end{equation}

\begin{rem}
We note that, for bounded and continuous functions $f$, Lemma~\ref{lemma:main} can also be derived from the
classical Crofton formula, well known in integral geometry (see for instance \cite{EisSul00}). In the appendix
we present an alternative proof of Lemma~\ref{lemma:main} constructed from Crofton formula.
\end{rem}

Let us now state our main result, which follows immediately from two applications of Lemma~\ref{lemma:main}.

\begin{thm}\label{thm:main}
Let $n\geqslant 2$ and let $f:\left]a,b\right[^{n+2}\to\mathbb{R}$ be an integrable function. Then we have
\begin{eqnarray*}
\lefteqn{\int_{\left]a,b\right[^n}f(\mathbf{x},\min x_i,\max x_i) \, d\mathbf{x}}\\
&=& \sum_{\textstyle{j,k=1\atop j\neq k}}^n \int_a^b dv\int_a^v du\int_{]u,v[^{n-2}} f(\mathbf{x},u,v\mid
x_j=u,x_k=v)\, \prod_{i\in [n]\setminus \{j,k\}} dx_i.
\end{eqnarray*}
\end{thm}

A direct use of this result leads to formula (\ref{eq:VarRan}). Indeed, as the integrand is symmetric in its
variables, we simply need to consider
$$
f(\mathbf{x},u,v\mid x_j=u,x_k=v) = \frac 1{v-u}\, s^2(x_1,\ldots,x_{n-2},u,v),
$$
where, for any fixed $a<u<v<b$, the right-hand side is a quadratic polynomial in $x_1,\ldots,x_{n-2}$.

\section{Application to aggregation function theory}

We now apply our main result to the computation of orness and andness average values of internal functions and
to the computation of idempotency average values of conjunctive and disjunctive functions.

\subsection{Internal functions}

We recall the concept of internal functions, which was introduced in the theory of means and aggregation
functions.

\begin{defn}
A function $F:\left]a,b\right[^n\to\mathbb{R}$ is said to be {\em internal}\/ if
$$
\min x_i \leqslant F(\mathbf{x}) \leqslant \max x_i \qquad (\mathbf{x}\in\left]a,b\right[^n).
$$
\end{defn}

Internality is a property introduced by Cauchy~\cite{Cau21} who considered in 1821 the {\em mean}\/ of $n$
independent variables $x_1,\ldots,x_n$ as a function $F(x_1,\ldots,x_n)$ which should be internal to the set of
$x_i$ values. Internal functions, also called Cauchy means, are very often encountered in the literature on
aggregation functions. Most of the classical means, such as the arithmetic mean, the geometric mean, and their
weighted versions, are Cauchy means. For overviews on means and aggregation functions, see the monograph
\cite{Bul03} and the edited book \cite{CalMayMes02}.

It is straightforward to see that a function $F:\left]a,b\right[^n\to\mathbb{R}$ is internal if and only if
there is a function $f$ from $\left]a,b\right[^n\setminus {\rm diag}(\left]a,b\right[^n)$ to $[0,1]$ such that
$$
F(\mathbf{x})=\min x_i+f(\mathbf{x})\, (\max x_i-\min x_i),
$$
where ${\rm diag}(\left]a,b\right[^n):=\{(x,\ldots,x)\in\left]a,b\right[^n : x\in\left]a,b\right[\}$.

Starting from this observation, Dujmovi\'c \cite{Duj73} (see also \cite{Duj05}) introduced the following
concepts of local orness and andness functions, rediscovered independently by Fern\'andez Salido and Murakami
\cite{FerMur03} as orness and andness distribution functions.

\begin{defn}
The {\em orness distribution function}\/ (resp.\ {\em andness distribution function}) associated with an
internal function $F:\left]a,b\right[^n\to\mathbb{R}$ is a function $odf_F$ (resp.~$adf_F$), from
$\left]a,b\right[^n\setminus {\rm diag}(\left]a,b\right[^n)$ to $\left[0,1\right]$, defined as
$$
odf_F(\mathbf{x})=\frac{F(\mathbf{x})-\min x_i}{\max x_i-\min x_i} \qquad
\mbox{(resp.~$\displaystyle{adf_F(\mathbf{x})=\frac{\max x_i-F(\mathbf{x})}{\max x_i-\min x_i}}$)}.
$$
\end{defn}

Thus defined, the orness distribution function (resp.\ andness distribution function) associated with an
internal function $F:\left]a,b\right[^n\to\mathbb{R}$ measures, at each $\mathbf{x}\in\left]a,b\right[^n$, the
extent to which $F(\mathbf{x})$ is close to $\max x_i$ (resp.\ $\min x_i$), that is, the extent to which
$F(\mathbf{x})$ has a disjunctive (resp.\ conjunctive) or orlike (resp.\ andlike) behavior.

To measure the average orness or andness quality of an internal function over its domain, Dujmovi\'c
\cite{Duj73} also introduced the concepts of mean local orness and andness, later called orness and andness
average values by Fern\'andez Salido and Murakami \cite{FerMur03}.

\begin{defn}
The {\em orness average value}\/ (resp.\ {\em andness average value}) of an internal and integrable function
$F:\left]a,b\right[^n\to\mathbb{R}$ is defined as
$$
\overline{odf}_F=\frac 1{(b-a)^n}\int_{]a,b[^n}odf_F(\mathbf{x})\, d\mathbf{x} \qquad
\mbox{(resp.~$\displaystyle{\overline{adf}_F=\frac 1{(b-a)^n}\int_{]a,b[^n}adf_F(\mathbf{x})\, d\mathbf{x}}$)}.
$$
\end{defn}

As an immediate property, we note that
$$
odf_F(\mathbf{x})+adf_F(\mathbf{x})=1,
$$
which entails $\overline{odf}_F+\overline{adf}_F=1$. Thus, as expected, both $\overline{odf}_F$ and
$\overline{adf}_F$ render the same information and hence we can restrict ourselves to the computation of
$\overline{odf}_F$.

Even though the computation of $\overline{odf}_F$ remains very difficult in most of the cases,
Theorem~\ref{thm:main} enables us to rewrite this integral in a more practical form, namely
$$
\overline{odf}_F = \frac 1{(b-a)^n}\sum_{\textstyle{j,k=1\atop j\neq k}}^n \int_a^b dv\int_a^v
du\int_{]u,v[^{n-2}} \frac{F(\mathbf{x}\mid x_j=u,x_k=v)-u}{v-u}\, \prod_{i\in [n]\setminus \{j,k\}} dx_i.
$$

The following two examples demonstrate the power of this formula:

\begin{exmp}
Let us calculate the orness average value over $[0,1]^n$ of the geometric mean
$$
G^{(n)}(\mathbf{x})=\prod_{i=1}^n x_i^{1/n}.
$$
The case $n=2$ is straightforward. Using (\ref{eq:nMaxMin}) with $f(u,v)=\frac{\sqrt{uv}-u}{v-u}$, we obtain
$\overline{odf}_{G^{(2)}} = \ln 4 - 1$.

Assume now that $n\geqslant 3$. As the integrand is a symmetric function, we can simply consider
$$
G^{(n)}(\mathbf{x}\mid x_j=u,x_k=v) = G^{(n)}(x_1,\ldots,x_{n-2},u,v)
$$
and hence, we have
\begin{eqnarray*}
\lefteqn{\int_{]u,v[^{n-2}} G^{(n)}(\mathbf{x}\mid x_j=u,x_k=v)\, \prod_{i\in [n]\setminus \{j,k\}} dx_i}\\
&=& \Big(\frac n{n+1}\Big)^{n-2}\,\big(v^{1+1/n}-u^{1+1/n}\big)^{n-2}u^{1/n}v^{1/n}.
\end{eqnarray*}
Then, using the binomial theorem and observing that $\frac 1{v-u}=\frac 1v\sum_{i=0}^{\infty}(\frac uv)^i$, we
obtain
\begin{eqnarray*}
\lefteqn{\int_0^1 dv\int_0^v \, \frac{\big(v^{1+1/n}-u^{1+1/n}\big)^{n-2}u^{1/n}v^{1/n}}{v-u}\, du}\\
&=& \sum_{i=0}^{\infty}\sum_{k=0}^{n-2}{n-2\choose k}\frac{(-1)^k}{n i+(k+1)(n+1)}\\
&=& \sum_{i=0}^{\infty}\sum_{k=0}^{n-2}{n-2\choose k}(-1)^k \int_0^1 x^{n i+k n+k+n}\, dx\\
&=& \int_0^1\frac{x^n(1-x^{n+1})^{n-2}}{1-x^n}\, dx.
\end{eqnarray*}
Finally,
$$
\overline{odf}_{G^{(n)}}=n(n-1)\Big(\frac{n}{n+1}\Big)^{n-2}\int_0^1\frac{x^n(1-x^{n+1})^{n-2}}{1-x^n}\, dx -
\frac 1{n-2}.
$$
The values of $\overline{odf}_{G^{(n)}}$ for $n=2,3,4,5$ are $\ln4 -1$, $\frac{\sqrt{3}\pi}2-\frac{47}{20}$,
$\frac{96\ln 2}{25}-\frac{8837}{3850}$, $\frac{25\pi}{27}\sqrt{\frac 52(25-11\sqrt{5})}-\frac{2454487}{960336}$,
respectively.
\end{exmp}

\begin{exmp}\label{ex:Choquet}
Let us calculate the orness average value over $[0,1]^n$ of a function of the form
$$
C_a^{(n)}(\mathbf{x})=\sum_{S\subseteq [n]}a(S)\min_{i\in S}x_i,
$$
where the set function $a:2^{[n]}\to\mathbb{R}$ fulfills
$$
a(\varnothing)=0 \quad\mbox{and}\quad \sum_{S\subseteq [n]}a(S)=1
$$
and is chosen so that the function $C_a^{(n)}$ is nondecreasing in each variable. Such a function is known in
aggregation function theory as a {\em Lov\'asz extension}\/ or a {\em discrete Choquet integral}\/ (see for
instance \cite{GraMarRou00,Mar00}). As particular cases, we can consider any weighted mean $\sum_iw_ix_i$ and
any convex combination $\sum_iw_ix_{(i)}$ of order statistics.

The case $n=2$ is easy. We simply obtain $\overline{odf}_{C_a^{(2)}} = \frac 12\big(a(\{1\})+a(\{2\})\big)$.

Assume now that $n\geqslant 3$. For any $0\leqslant u<v\leqslant 1$ and any $\mathbf{x}\in [u,v]^{n-2}$, we have
$$
C_a^{(n)}(\mathbf{x}\mid x_j=u,x_k=v)=\sum_{S\ni j} a(S) u + \sum_{\textstyle{S\not\ni j\atop S\ni k}}
a(S)\min_{i\in S\setminus\{k\}} x_i + \sum_{\textstyle{S\not\ni j\atop S\not\ni k}} a(S)\min_{i\in S} x_i.
$$
Setting $s:=|S|$, from (\ref{eq:CaseMinSab}) it follows that
\begin{eqnarray*}
\lefteqn{\int_{]u,v[^{n-2}} C_a^{(n)}(\mathbf{x}\mid x_j=u,x_k=v)\, \prod_{i\in [n]\setminus \{j,k\}} dx_i}\\
&=& (v-u)^{n-2}\bigg(\sum_{S\ni j} a(S) u + \sum_{\textstyle{S\not\ni j\atop S\ni k}} a(S)\frac{u(s-1)+v}{s} +
\sum_{\textstyle{S\not\ni j\atop S\not\ni k}} a(S)\frac{u s+v}{s+1}\bigg).
\end{eqnarray*}
Then, we obtain
\begin{eqnarray*}
\lefteqn{\int_0^1 dv \int_0^v \frac{du}{v-u}\int_{]u,v[^{n-2}} C_a^{(n)}(\mathbf{x}\mid x_j=u,x_k=v)\, \prod_{i\in [n]\setminus \{j,k\}} dx_i}\\
&=& \frac 1{n(n-1)(n-2)} \bigg(\sum_{S\ni j} a(S) + \sum_{\textstyle{S\not\ni j\atop S\ni k}} a(S)\frac{n+s-2}s
+ \sum_{\textstyle{S\not\ni j\atop S\not\ni k}} a(S)\frac{n+s-1}{s+1}\bigg).
\end{eqnarray*}
Summing over $j,k=1\ldots,n$, with $j\neq k$, and then rearranging the terms we finally obtain
\begin{eqnarray*}
\overline{odf}_{C_a^{(n)}}
&=& \biggl(\sum_{S\subseteq [n]} a(S)\,\frac{n(n-1)+s-1}{(n-1)(n-2)(s+1)}\biggr)-\frac 1{n-2}\\
&=& \sum_{S\subseteq [n]} a(S)\,\biggl(\frac{n(n-1)+s-1}{(n-1)(n-2)(s+1)}-\frac 1{n-2}\biggr)\\
&=& \frac 1{n-1}\sum_{S\subseteq [n]} a(S)\,\frac{n-s}{s+1},
\end{eqnarray*}
which includes the case $n=2$.
\end{exmp}

To overcome the difficulty of calculating intractable orness average values, Dujmovi\'c \cite{Duj74} introduced
the next concept of global orness and andness measures (see also \cite{Duj05,FerMur03}). Denote by
$\overline{F}$ the average value of any internal and integrable function $F:\left]a,b\right[^n\to\mathbb{R}$
over its domain, that is,
$$
\overline{F}:=\frac 1{(b-a)^n}\int_{]a,b[^n}F(\mathbf{x})\, d\mathbf{x}.
$$
\begin{defn}
The {\em global orness value}\/ (resp.\ {\em global andness value}) of an internal and integrable function
$F:\left]a,b\right[^n\to\mathbb{R}$ is defined as
$$
{\rm orness}_F=\frac{\overline{F}-\overline{{\rm Min}}}{\overline{{\rm Max}}-\overline{{\rm Min}}} \qquad
\mbox{(resp.~$\displaystyle{{\rm andness}_F=\frac{\overline{{\rm Max}}-\overline{F}}{\overline{{\rm
Max}}-\overline{{\rm Min}}}}$)},
$$
where ${\rm Min}$ and ${\rm Max}$ are, respectively, the minimum and maximum functions defined in
$\left]a,b\right[^n$.
\end{defn}

For example, considering the geometric mean $G^{(n)}(\mathbf{x})=\prod_{i=1}^n x_i^{1/n}$ in $[0,1]^n$, we
simply obtain
$$
{\rm orness}_{G^{(n)}} = -\frac 1{n-1}+\frac{n+1}{n-1}\,\overline{G^{(n)}} = -\frac
1{n-1}+\frac{n+1}{n-1}\Big(\frac n{n+1}\Big)^n.
$$
Considering the discrete Choquet integral $C_a^{(n)}$ in $[0,1]^n$, as defined in Example~\ref{ex:Choquet}, we
get
\begin{eqnarray*}
{\rm orness}_{C_a^{(n)}} &=& -\frac 1{n-1}+\frac{n+1}{n-1}\,\overline{C_a^{(n)}} =  -\frac 1{n-1}+\frac{n+1}{n-1}\,\sum_{S\subseteq [n]} a(S)\, \frac 1{|S|+1}\\
&=& \frac 1{n-1}\sum_{S\subseteq [n]} a(S)\Big(\frac{n+1}{|S|+1}-1\Big)\\
&=& \frac 1{n-1}\sum_{S\subseteq [n]} a(S)\,\frac{n-|S|}{|S|+1}.
\end{eqnarray*}

Surprisingly enough, in $[0,1]^n$ we have
$$
{\rm orness}_{C_a^{(n)}}=\overline{odf}_{C_a^{(n)}},
$$
that is, for any discrete Choquet integral, the global orness value identifies with the orness average value, a
result already reached by Fern\'andez Salido and Murakami \cite{FerMur03} for the special case of symmetric
Choquet integrals, that is, convex combinations of order statistics.

The interesting question of determining those internal functions $F:\left]a,b\right[^n\to\mathbb{R}$ fulfilling
the equation ${\rm orness}_{F}=\overline{odf}_{F}$ remains open.

\subsection{Conjunctive and disjunctive functions}

Let us now consider conjunctive and disjunctive functions.

\begin{defn}
A function $F:\left]a,b\right[^n\to\R$ is said to be {\em conjunctive}\/ (resp.\ {\em disjunctive}) if
$$
a\leqslant F(\mathbf{x})\leqslant\min x_i \qquad \big(\mbox{resp.\ }\max x_i\leqslant F(\mathbf{x})\leqslant
b\big).
$$
\end{defn}

Prominent examples of conjunctive (resp.\ disjunctive) functions in the literature are {\em t-norms}\/ (resp.\
{\em t-conorms}), which are symmetric, associative, and nondecreasing functions, from $[0,1]^2$ to $[0,1]$, with
$0$ (resp.\ $1$) as the neutral element. For an account on t-norms and t-conorms, see for instance the book by
Alsina et al.\ \cite{AlsFraSch06}.

Clearly, a function $F:\left]a,b\right[^n\to\R$ is conjunctive (resp.\ disjunctive) if and only if there is a
function $f:\left]a,b\right[^n\to [0,1]$ such that
$$
F(\mathbf{x})=a+f(\mathbf{x})(\min x_i-a) \qquad \big(\mbox{resp.\ }F(\mathbf{x})=b-f(\mathbf{x})(b-\max
x_i)\big).
$$

Just as for the orness and andness distribution functions, we can naturally define the concept of idempotency
distribution function associated with a conjunctive (resp.\ disjunctive) function $F:\left]a,b\right[^n\to\R$ as
a measure, at each $\mathbf{x}\in\left]a,b\right[^n$, of the extent to which $F$ is idempotent (i.e., such that
$F(x,\ldots,x)=x$), that is, the extent to which $F$ is close to $\min x_i$ (resp.\ $\max x_i$).

\begin{defn}
The {\em idempotency distribution function}\/ associated with a conjunctive (resp.\ disjunctive) function
$F:\left]a,b\right[^n\to\R$ is a function $idf_F:\left]a,b\right[^n\to [0,1]$, defined as
$$
idf_F(\mathbf{x})=\frac{F(\mathbf{x})-a}{\min x_i-a} \qquad
\mbox{(resp.~$\displaystyle{idf_F(\mathbf{x})=\frac{b-F(\mathbf{x})}{b-\max x_i}}$)}.
$$
\end{defn}

We can now introduce the concept of idempotency average value as follows.

\begin{defn}
The {\em idempotency average value}\/ of a conjunctive or disjunctive function $F:\left]a,b\right[^n\to\R$ is
defined as
$$
\overline{idf}_F = \frac 1{(b-a)^n}\,\int_{]a,b[^n} idf_F(\mathbf{x})\, d\mathbf{x}.
$$
\end{defn}

According to Lemma~\ref{lemma:main}, for any conjunctive function $F:\left]a,b\right[^n\to\R$ for instance, we
can write
$$
\overline{idf}_F = \frac 1{(b-a)^n}\, \sum_{j=1}^n \int_a^b du \int_{]u,b[^{n-1}} \frac{F(\mathbf{x}\mid
x_j=u)-a}{u-a}\,\prod_{i\in [n]\setminus\{j\}} dx_i.
$$

The following concept of global idempotency value was introduced by Koles\'arov\'a \cite{Kol06} for t-norms as
an idempotency measure:

\begin{defn}
The {\em global idempotency value}\/ of a conjunctive (resp.\ disjunctive) function $F:\left]a,b\right[^n\to\R$
is defined by
$$
{\rm idemp}_F = \frac{\overline{F}-a}{\overline{{\rm Min}}-a} \qquad \mbox{(resp.~$\displaystyle{{\rm idemp}_F =
\frac{b-\overline{F}}{b-\overline{{\rm Max}}}}$)}.
$$
\end{defn}

\begin{exmp}
Let us calculate the idempotency average value and the global idempotency value over $[0,1]^n$ of the product
$$
P^{(n)}(\mathbf{x})=\prod_{i=1}^n x_i,
$$
which is a conjunctive function.

We immediately obtain
\begin{eqnarray*}
\overline{idf}_{P^{(n)}} &=& n\int_0^1 du \int_{[u,1]^{n-1}}\Big(\prod_{i=1}^{n-1} x_i\Big)\, dx_1\cdots dx_{n-1}\\
&=& \frac{n}{2^{n-1}} \int_0^1 (1-u^2)^{n-1}\, du.
\end{eqnarray*}
Setting $u=v^{1/2}$ and then using the classical beta function $$\mathrm{B}(a,b)=\int_0^1 t^{a-1}(1-t)^{b-1}\,
dt,$$ we obtain

\begin{eqnarray*}
\overline{idf}_{P^{(n)}} &=& \frac{n}{2^{n}} \int_0^1 v^{-1/2}(1-v)^{n-1}\, dv \\
&=& \frac{n}{2^{n}} \ \mathrm{B}(1/2,n) = \frac{n}{2^{n}}\, \frac{\Gamma(1/2)\Gamma(n)}{\Gamma(n+1/2)}\\
&=& \frac{2^{n-1}}{{2n-1 \choose n}}.
\end{eqnarray*}

On the other hand, we have
$$
{\rm idemp}_{P^{(n)}} = (n+1)\int_{[0,1]^{n}} \prod_{i=1}^n x_i\, d\mathbf{x} = \frac{n+1}{2^n}.
$$
\end{exmp}

\section{Application to probability theory}

Since the idea behind our results merely consists in breaking the integration domain into smaller regions,
Lemma~\ref{lemma:main} and Theorem~\ref{thm:main} can be straightforwardly extended to Riemann-Stieltjes
integrals, thus making it possible to consider average values from various probability distributions.

Consider a measurable function $g:\R^{n+2}\to\R$ and $n$ independent random variables $X_1,\ldots,X_n$, $X_i$
$(i\in [n])$ having distribution function $F_i(x)$. Define the random variable $Y_g$ as
$$
Y_g:=g(\mathbf{X},\min X_i,\max X_i),
$$
where $\mathbf{X}$ denotes the vector $(X_1,\ldots,X_n)$.

The direct generalization of Theorem~\ref{thm:main} to Riemann-Stieltjes integrals can be used to evaluate the
expected value of $Y_g$, namely
$$
\mathbf{E}[Y_g]=\int_{\R^n}g(\mathbf{x},\min x_i,\max x_i)\, dF_1(x_1)\cdots dF_n(x_n).
$$
It can also be used in the evaluation of the distribution function of $Y_g$, which is defined as
\begin{eqnarray*}
F_g(z) &=& \mathbf{E}[H(z-Y_g)]\\
&=& \int_{\R^n}H\big(z-g(\mathbf{x},\min x_i,\max x_i)\big)\, dF_1(x_1)\cdots dF_n(x_n),
\end{eqnarray*}
where $H:\R\to\{0,1\}$ is the Heaviside step function, defined by $H(x)=1$ if $x\geqslant 0$, and 0 otherwise.
Note that the case where $Y_g$ is a lattice polynomial (max-min combination) of the variables $X_1,\ldots,X_n$
has been thoroughly investigated by the author in \cite{Mar06,Mard}.

To keep our exposition simple, let us examine the special case where $Y_g$ depends only on $\min X_i$ and $\max
X_i$, that is,
$$
Y_g:=g(\min X_i,\max X_i),
$$
where $g:\R^2\to\R$ is a measurable function. In this case, our method immediately leads to
\begin{eqnarray*}
\mathbf{E}[Y_g] &=& \sum_{\textstyle{j,k=1\atop j\neq k}}^n \int_{-\infty}^{\infty} dF_k(v)\int_{-\infty}^v
g(u,v)\, dF_j(u)\int_{]u,v[^{n-2}}
\prod_{i\in [n]\setminus \{j,k\}} dF_i(x_i)\\
&=& \sum_{\textstyle{j,k=1\atop j\neq k}}^n \int_{-\infty}^{\infty} dF_k(v)\int_{-\infty}^v g(u,v)\, \prod_{i\in
[n]\setminus \{j,k\}}\big(F_i(v)-F_i(u)\big) \, dF_j(u).
\end{eqnarray*}

In the particular case where the random variables $X_1,\ldots,X_n$ are independent and identically distributed,
each with distribution function $F(x)$, the expected value clearly reduces to
\begin{equation}\label{eq:Ygiid}
\mathbf{E}[Y_g] = n(n-1) \int_{-\infty}^{\infty} dF(v)\int_{-\infty}^v g(u,v)\, \big(F(v)-F(u)\big)^{n-2} \,
dF(u),
\end{equation}
which generalizes formula (\ref{eq:nMaxMin}).

\begin{exmp}
For exponential variables $X_1,\ldots,X_n$, each with distribution function $F(x)=1-e^{-\lambda x}$ $(x>0)$, we
simply have
$$
\mathbf{E}[Y_g] = n(n-1) \int_{0}^{\infty} \lambda\, e^{-\lambda v}\, dv \int_{0}^v g(u,v)\, \big(e^{-\lambda
u}-e^{-\lambda v}\big)^{n-2} \, \lambda\, e^{-\lambda u}\, du.
$$
Using the change of variables $x= e^{-\lambda u}$ and $y = e^{-\lambda u}-e^{-\lambda v}$, this integral can be
easily rewritten as
$$
\mathbf{E}[Y_g] = n(n-1) \int_{0}^{1}  dx \int_{0}^{x} y^{n-2}\, g\big(-\frac 1{\lambda}\,\ln(x),-\frac
1{\lambda}\,\ln(x-y)\big)\, dy.
$$
\end{exmp}

\begin{exmp}
Let us calculate the distribution function and the raw moments of the random variable
$$
Y=\frac{\max X_i-\min X_i}{\max X_i}
$$
from the uniform distribution over $]0,1]^n$.

The raw moments can be calculated very easily from (\ref{eq:Ygiid}). For any integer $r\geqslant 0$, we have
$$
\mathbf{E}[Y^r]=n(n-1) \int_{0}^{1} dv\int_{0}^v \Big(\frac{v-u}v\Big)^r\, (v-u)^{n-2} \, du=\frac{n-1}{n+r-1}.
$$
On the other hand, the distribution of $Y$ is simply given by
$$
F(z)= n(n-1) \int_{0}^{1} dv\int_{0}^v H\big(z-\frac{v-u}v\big)\, (v-u)^{n-2} \, du,
$$
that is,
$$
F(z)=
\begin{cases}
0, & \mbox{if $z\leqslant 0$,}\\
\displaystyle{n(n-1) \int_{0}^{1} dv\int_{v(1-z)}^v (v-u)^{n-2} \, du} = z^{n-1}, & \mbox{if $0\leqslant z\leqslant 1$,}\\
\displaystyle{n(n-1) \int_{0}^{1} dv\int_{0}^v (v-u)^{n-2} \, du = 1}, & \mbox{if $1\leqslant z$.}
\end{cases}
$$
\end{exmp}

\appendix
\section*{Appendix: The use of Crofton formula}

We provide a proof of Lemma~\ref{lemma:main} as a direct consequence of Crofton formula. See \cite{EisSul00} for
a very good expository note on Crofton formula.

For $0<v<v+h<V$, let $D(v)$ be a domain of area or volume $v$. By a domain we mean a closed bounded convex set
in $\R^k$ for some $k$. Assume that for $v_1<v_2$ we have $D(v_1)\subset D(v_2)$. Let $X_1,\ldots,X_n$ be $n$
independent points randomly selected with uniform distribution in $D(v+h)$ and let $Y=f(X_1,\ldots,X_n)$, where
$f$ is a bounded function. Let $A(v)$ be the event that all the points are in $D(v)$ and let $B_j(v,h)$ be the
event that $X_j\in D(v+h)-D(v)$ and $X_i\in D(v)$ for all $i\neq j$. Let $\mu(v)=\mathbf{E}[Y|A(v)]$ and let
$\mu_j^*(v,h)=\mathbf{E}[Y|B_j(v,h)]$.

In its nonsymmetric version, Crofton formula states that, if $\lim_{h\to 0}\mu_j^*(v,h)=\mu_j(v)$ exists and is
continuous for all $j$, then for $V>0$,
$$
\mathbf{E}[Y]=\mu(V)=\frac 1{V^n}\sum_{j=1}^n\int_0^Vv^{n-1}\mu_j(v)\, dv.
$$

Choosing $V=b-a$ and $D(v)=[a,a+v]$ (which implies $D(V)=[a,b]$), we simply obtain
\begin{eqnarray*}
\mu(V)&=&\frac 1{(b-a)^n}\int_{[a,b]^n}f(\mathbf{x})\, d\mathbf{x},\\
\mu_j^*(v,h) &=& \frac 1{v^{n-1}h}\int_a^{a+v}\, dx_1\cdots\int_{a+v}^{a+v+h}\, dx_j\cdots\int_a^{a+v}f(\mathbf{x})\, dx_n,\\
\mu_j(v) &=& \frac 1{v^{n-1}}\int_{[a,a+v]^{n-1}}f(\mathbf{x}\mid x_j=a+v)\, \prod_{i\in [n]\setminus \{j\}}
dx_i.
\end{eqnarray*}
If $f$ is continuous in each argument then $\lim_{h\to 0}\mu_j^*(v,h)=\mu_j(v)$ exists and is continuous.
According to Crofton formula, we obtain
$$
\mu(V)=\frac 1{V^n}\sum_{j=1}^n\int_0^V\, dv\, \int_{[a,a+v]^{n-1}}f(\mathbf{x}\mid x_j=a+v)\, \prod_{i\in
[n]\setminus \{j\}} dx_i,
$$
which proves the second formula of Lemma~\ref{lemma:main}. The first one can be established similarly by
considering $D(v)=[b-v,b]$. \qed



\end{document}